\documentstyle{article}
\newcommand{\PP}{{\bf P}}
\newcommand{\PN}{\PP^N}
\renewcommand{\O}{{\cal O}}
\renewcommand{\S}{{\bf S}}
\newcommand{\Proj}{\mathop{\rm Proj}\nolimits}
\newcommand{\codim}{\mathop{\rm codim}\nolimits}
\newcommand{\idsheaf}[1]{{\cal I}_{#1}}
\newtheorem{prop}{Proposition}[section]
\newtheorem{lemma}[prop]{Lemma}
\newcommand{\Proof}{\par\addvspace{\smallskipamount}
                    \noindent{\bf Proof }}
\begin{document}
\author{S.~L'vovsky\thanks{Partially supported by an ISF grant}}
\title{On Landsberg's criterion for complete intersections}
\date{August 24, 1994}
\maketitle
\section*{Introduction}
In his preprint~\cite{Lan}, J.M.~Landsberg introduces an
elementary characterization of complete intersections
(Proposition~1.2 in \cite{Lan}). The proof of this proposition
uses the method of moving frames. The aim of this note is to
present an elementary proof of Landsberg's criterion that is
valid over any ground field.

\section{Notation and statement of results}
Let $k$ be an algebraically closed field and $\PP^N= \Proj
k[T_0,\ldots,T_N]$ the $N$-dimensional projective space over $k$.
If $F$ is a homogeneous polynomial in $T_0,\ldots ,T_N$, we will
denote by $Z(F)\subset \PP^N$ the hypersurface defined by $F$. If
$F$ is a homogeneous polynomial and $x=(x_0:\ldots:x_N)\in \PN$,
put $d_x F=\left(\partial F/\partial T_0(z),\ldots, \partial
F/\partial T_N(z)\right)\in k^{N+1}$ (actually $d_x F$ depends on
the choice of homogeneous coordinates for $x$; this abuse of
notation should not lead to confusion). If $x\in X$, where
$X\subset\PN$ is a projective variety, then $T_xX\subset \PN$
denotes the embedded Zariski tangent space to $X$ at $x$.

If $X\subset \PP^N$ is a projective variety, then its ideal sheaf
will be denoted by $\idsheaf X\subset \O_{\PP^N}$ and its
homogeneous ideal by $I_X\subset k[T_0, \ldots, T_N]$. We will
say that a hypersurface $Y=Z(F)$ {\em trivially contains $X$\/}
iff $F=\sum G_iF_i$, where $G_i$'s and $F_i$'s are homogeneous
polynomials, $F_i$ vanish on $X$ for all $i$, and $\deg F_i<\deg
F$ for all $i$. If $Y$ trivially contains $X$, then $Y\supset X$.
We will say that a hypersurface $W$ {\em non-trivially contains
$X$\/} iff $W$ contains $X$, but not trivially.

The following proposition is a slight reformulation of
Landsberg's criterion (cf.\ \cite[Proposition 1.2]{Lan}):
\begin{prop}
For a projective variety $X\subset \PP^N$, the following
conditions are equivalent:
\begin{itemize}
\item[(i)]
$X$ is a complete intersection.
\item[(ii)]
There exists a smooth point $x\in X$ having the following
property: any hypersurface $W\subset \PP^N$ that non-trivially
contains $X$ must be smooth at $x$.
\item[(iii)]
For any smooth point $x\in X$ and any hypersurface $W$ that
non-trivially contains $X$, $W$ is smooth at $x$.
\item[(iv)]
For any smooth point $x\in X$ and any hypersurface $W$ that
non-trivially contains $X$, $T_xW$ cannot contain an intersection
$\bigcap_i T_xW_i$, where each $W_i$ is a hypersurface s.t.\ $W_i
\supset X$ and $\deg W_i<\deg W$ (it is understood that the
intersection of an empty family of tangent spaces is the entire
$\PN$).
\end{itemize}
\end{prop}

\section{Proofs}
For the sequel we need two lemmas.
\begin{lemma}\label{subst}
Let $F_1,\ldots,F_r$ be homogeneous polynomials over $k$ in
$T_0,\ldots, T_N$. Assume that $x=(x_0:\ldots:x_N)\in \PN$ is
their common zero and that the vectors $d_xF_1,\ldots, d_xF_r$
are linearly dependent. Then one of the following alternatives
holds:
\begin{enumerate}
\item
There is $j\in [1;r]$ s.t.\ $F_j$ belongs to the ideal in
$k[T_0,\ldots,T_N]$ generated by $F_i$'s with $i\ne j$.
\item
There are homogeneous polynomials $\tilde F_0,\ldots, \tilde F_N$
s.t.\ the ideals $(F_0,\ldots,F_N)$ and $(\tilde F_0,\ldots,
\tilde F_N)$ coincide, $\deg \tilde F_i=\deg F_i$ for all $i$,
and $d_x \tilde F_j=0$ for some $j$.
\end{enumerate}
\end{lemma}
{\bf Proof.}
Let the shortest linear relation among $d_xF_j$'s have the form
$$
  \lambda_1d_xF_1+\cdots+\lambda_sd_xF_s=0,
$$
where $\lambda_j\ne 0$ for all $j$. Reordering $F_j$'s if
necessary, we may assume that $\deg F_1\le \deg F_2\le \cdots\le
\deg F_s$. Let $t$ be such a number that $\deg F_t=\deg F_s$ and
$\deg F_{t-1} <\deg F_s$ (if $\deg F_1=\deg F_s$, set $t=1$).

If the polynomials $F_t,\ldots,F_s$ are linearly dependent, then
it is clear that one of them lies in the ideal generated by the
others and there is nothing more to prove. Assume from now on
that $F_t, F_{t+1},\ldots, F_s$ are linearly independent. Then
there exists an index $j\in[t;s]$ and numbers $\mu_i$, where $i
\in [t;s]$ s.t.\
\begin{equation}\label{G:def}
  F_j=\sum_{i\in [t;s]\setminus \{j\}}\mu_i F_i+
        \mu_j(\lambda_t F_t+\cdots+\lambda_s F_s).
\end{equation}
For each $i\in[1;t-1]$, choose a homogeneous polynomial $G_i$
s.t.\ $\deg G_i= \deg F_s-\deg F_i$ and $G_i(x_0,\ldots,
x_N)=\lambda_i$, and set
\begin{equation}\label{tilde:def}
  \tilde F_j=\sum_{i<t}G_iF_i+\sum_{i\ge t}\lambda_i F_i.
\end{equation}
If $\tilde F_j=0$, then $F_s\in (F_1,\ldots, F_{s-1})$ and the
first alternative holds. Otherwise, $\deg \tilde F_j=\deg F_j$,
$d_x\tilde F_j=0$ by virtue of (\ref{tilde:def}), and it
follows from (\ref{G:def}) and (\ref{tilde:def}) that
$$
  F_j=\sum_{i\in [t;s]\setminus \{j\}}\mu_i F_i +\mu_j \tilde F_j
  -\mu_j\sum_{i<t}G_i F_i,
$$
whence $(F_1,\ldots,F_{j-1}, \tilde F_j, F_{j+1},\ldots,
F_s)=(F_1,\ldots,F_s)$. Hence in this case the second alternative
holds, and we are done.

The second lemma belongs to folklore. To state this lemma, let us
introduce some notation. Denote by $\S$ the set of sequences of
non-negative integers $\delta=(\delta_1,\delta_2,\ldots)$ s.t.\
$\delta_M=0$ for all $M\gg 0$. If $\delta,\eta\in \S$, we will
write~$\delta \succ \eta$ iff there is an integer $i$ s.t.\
$\delta_i >\eta_i$ and $\delta_j=\eta_j$ for all $j>i$.
\begin{lemma}\label{folk}
Any sequence $\delta_1 \succ \delta_2 \succ\cdots$ must
terminate.
\end{lemma}
\noindent {\bf Proof.}
For any $\delta\in \S$, set $n(\delta)=\max\{j:\delta_j\ne 0\}$,
$\ell(\delta)=\delta_{n(\delta)}> 0$. If $\delta\succ \eta$ and
$n(\delta)=n(\eta)$, then $\ell(\delta)\ge \ell(\eta)$. Let us
prove the lemma by induction on $n(\delta_1)$.

If $n(\delta_1)\le 1$, the result is evident. Assuming that the
lemma is true whenever $n(\delta_1)< m$, suppose that there is an
infinite sequence $\delta_1 \succ \delta_2 \succ\cdots$ with
$n(\delta_1)=m$. If $n(\delta_j)<n(\delta_1)$ for some $j$, we
arrive at a contradiction by the induction hypothesis. Hence,
$n(\delta_j) =n(\delta_1) =m$ for all $j$ and $\ell(\delta_1)\ge
\ell(\delta_2)\ge\cdots >0$. Thus there exists an integer $N$
s.t\ $\ell(\delta_j)$ is connstant for $j\ge N$. For any $j\ge
N$, denote by $\delta'_j\in \S$ a sequence that is obtained from
$\delta_j$ by replacing its last positive term by zero. It is
clear that $\delta'_N \succ \delta'_{N+1} \succ\cdots$, and this
sequence is infinite by our assumption. This is again impossible
by the induction hypothesis since $n(\delta'_j)< n(\delta_j)=m$,
whence the lemma.

\smallskip
\Proof of $(ii)\Rightarrow (i)$.
Put $a=N-\dim X$. Let $(F_1,\ldots,F_r)$ be a system of
(homogeneous) generators of $I_X$. To any such system assign a
sequence $\delta(F_1,\ldots, F_r) \in \S$, where $\delta(F_1,
\ldots,F_r)_i = \#\{j\in [1;r]:\deg F_j=i\}$. I claim that
\begin{quote}
if $r>a$, then $I_X=(\Phi_1,\ldots,\Phi_s)$, where $\Phi_i$'s are
homogeneous polynomials s.t.\ $\delta(F_1,\ldots,F_r)\succ
\delta(\Phi_1,\ldots, \Phi_s)$.
\end{quote}
To prove this claim, observe that $d_xF_1,\ldots,d_xF_r$ are
linearly dependent since $X$ is smooth at $x$ and $r>\codim X$.
Now Lemma~\ref{subst} implies that either one of the $F_j$'s
(say, $F_1$) can be removed without affecting $I_X$, or
$I_X=(\tilde F_1,\ldots,\tilde F_r)$, where $\deg \tilde F_j=\deg
F_j$ for all $j$ and $d_x \tilde F_j=0$ for some $j$. In the
first case, the required $\Phi_1,\ldots,\Phi_s$ can be obtained
by merely removing $F_1$; in the second case, hypothesis~$(ii)$
shows that $\tilde F_j=\sum_{i=1}^t G_i\Psi_i$, where $\Psi_i\in
I_X$ and $\deg\Psi_i < \deg \tilde F_j$ for all $j$. Replacing
$\tilde F_j$ by $\Psi_1,\ldots, \Psi_t$ in the sequence $\tilde
F_1,\ldots,\tilde F_r$ and putting $s=r+t-1$, we obtain a
sequence $\Phi_1,\ldots,\Phi_s$ s.t.\ $I_X =(\Phi_1,\ldots,
\Phi_s)$ and $\delta(\tilde F_1,\ldots,\tilde F_r)\succ
\delta(\Phi_1,\ldots, \Phi_s)$. Since the degrees of $\tilde
F_j$'s and $F_j$'s are the same, this means that $\delta(F_1,
\ldots,F_r)\succ \delta(\Phi_1,\ldots, \Phi_s)$ as well, and the
claim is proved.

Now we can finish the proof as follows. If $r=a$, then $X$ is the
complete intersection of $Z(F_1),\ldots,Z(F_r)$ and there is
nothing to prove. If $r>a$, then by virtue of our claim we can
replace the system of generators $F_1,\ldots,F_r$ by
$\Phi_1,\ldots, \Phi_s$. Let us iterate this process. By virtue
of Lemma~\ref{folk} this process must terminate and by virtue of
our claim this is possible only when we have found a system of
exactly $a$ generators of the ideal $I_X$. This means that $X$ is
a complete intersection, thus completing our proof.

\smallskip

\Proof of $(iv)\Rightarrow
(iii)\Rightarrow (ii)$. Trivial.

\Proof of $(i)\Rightarrow (iv)$.
Let $X$ be a complete intersection of the hypersurfaces $Z(F_1),
\ldots, Z(F_a)$. Assume that a hypersurface $W=Z(F)$, with $F$
irreducible, non-trivially contains $X$ and that $x=(x_0:\ldots
:x_N)\in \PN$ is a smooth point of $X$; set $m=\deg F$. Since
$Z(F)\supset X$ and $X$ is a complete intersection of the
$Z(F_i)$'s, we see that
\begin{equation}\label{expr}
F=\sum G_iF_i;
\end{equation}
since $W$ contains $X$ non-trivially, at least some of the
$G_j$'s must be non-zero constants. Reordering $F_j$'s if
necessary, we may assume that $G_j$ is a constant (hence, $\deg
F_j=m$) iff $1\le j\le s$. Taking $d_x$ of the both parts of
(\ref{expr}), we see that
\begin{equation}\label{diffls}
d_xF=\sum_{i=1}^a c_i d_x F_i,\qquad
\mbox{where $c_i\ne 0$ for some $i\in [1;s]$.}
\end{equation}

On the other hand, assume that $W_i=Z(B_i)$ with irreducible
$B_i$'s. Then the hypothesis implies that $d_x F$ is a linear
combination of $d_xB_j$'s, and the fact that $X$ is a complete
intersection of $Z(F_t)$'s and $Z(B_j)\supset X$ implies that,
for each $j$, there is a relation
\begin{equation}\label{expr'}
B_j=\sum_{t>s} G_{jt}F_t
\end{equation}
(it suffices to sum only over $t>s$ since for $t\le s$ we have
$\deg F_t= \deg W > \deg B_j$). If we take $d_x$ of both parts of
(\ref{expr'}), we see that, for each $j$, $d_xB_j$ is a linear
combination of $d_x F_t$'s with $t>s$. Hence $d_xF$ is also a
linear combination of $d_x F_t$'s with $t>s$. Taking into account
(\ref{diffls}) we see that $d_xF_i$'s are linearly dependent.
This is, however, impossible since $x$ is a smooth point of the
comlete intersection of $Z(F_j)$'s. This contradiction completes
the proof.

\begin{tabbing}
Email address: \=serge@schcorr.msk.su\\
\>nskcsmoscow@glas.apc.org
\end{tabbing}

\begin{thebibliography}{1}
\bibitem{Lan}
Landsberg, J.~M. Differential-Geometric Characterizations of
Complete Intersections. Preprint, alg-geom/9407002
\end{thebibliography}
\end{document}